\documentclass[12pt]{amsart}
\usepackage{amssymb}
\usepackage{mathrsfs}

\theoremstyle{plain}

\newtheorem{theorem}{Theorem}
\newtheorem{construction}{Construction}
\newtheorem{lemmaX}{Lemma}

\newtheorem{example}{Example}
\newtheorem{remark}{Remark}
\theoremstyle{definition}

\newcommand{\eps}{\varepsilon}

\newcommand{\dist}{{\rm dist\,}}

\newcommand{\co}{{\rm conv\,}}

\newcommand{\spn}{{\rm span\,}}

\newcommand{\Aa}{\mathcal A}

\newcommand{\DD}{\mathcal D}

\newcommand{\C}{\mathbb C}
\newcommand{\R}{\mathbb R}

\newcommand{\NN}{\mathbb N}

\newcommand{\SSS}{\mathcal S}

\newif\ifComplain
\Complainfalse          %  suppress the \complain{blabla} remarks
\def\complain#1{\ifComplain\ifhmode \newline\fi{\sf *** \ \ #1
\\}\fi}

\newif\ifmarglab
\makeatletter
\def\label#1{\@bsphack\ifmarglab\marginpar{LAB:#1}\fi\if@filesw {\let\thepage\relax
   \def\protect{\noexpand\noexpand\noexpand}%
   \edef\@tempa{\write\@auxout{\string
      \newlabel{#1}{{\@currentlabel}{\thepage}}}}%
   \expandafter}\@tempa
   \if@nobreak \ifvmode\nobreak\fi\fi\fi\@esphack}
\makeatother
\marglabfalse
\usepackage{graphicx}
\long\def\onefigure#1#2{%  #1 picture,  #2  caption
\begin{figure*}[tbp]
\begin{center}
#1
\end{center}
\caption{#2}
\end{figure*}
} %end onefigure def

\newcommand\newipefig[2]%  #1 file name & label, # 2 caption
{\onefigure{\includegraphics{#1.pdf}}{\label{#1} #2}}
\begin{document}

\title[Sequences of $ m$-term deviations in Hilbert space]{Sequences of $m$-term deviations\\ in Hilbert space}

\author[P. A. Borodin]{Petr A. Borodin}
\address{Department of Mechanics and
Mathematics, Moscow State University, Moscow 119991, Russia}
\address{and}
\address{Moscow Center for Fundamental and Applied Mathematics}
\email{pborodin@inbox.ru}
\author[E. Kopeck\'a]{ Eva Kopeck\'a}
\address{Department of Mathematics\\
   University of Innsbruck\\
 A-6020 Innsbruck, Austria}
 \address{and}
\address{Moscow Center for Fundamental and Applied Mathematics}
\email {eva.kopecka@uibk.ac.at}

\thanks{The first author was supported by the grant of the Government of the
Russian Federation (project 14.W03.31.0031).}


\date{\today}

\begin{abstract}
Let $D$ be a dictionary in a Hilbert space $H$, that is, a set of unit elements whose linear combinations are dense in $H$. We consider the least $m$-term deviation $\sigma_m(x)$ of an element $x\in H$: this is  the distance of $x$ from the set of all   $m$-term linear combinations of elements of $D$. We prove a dichotomy result: for any dictionary $D$, either the sequence  $\{\sigma_m(x)\}_{m=0}^{\infty}$ decreases exponentially for every $x\in H$, or the rate of convergence $\sigma_m(x)\to 0$ can be arbitrarily slow. We seek universal dictionaries realizing all strictly decreasing  null sequences as sequences of $m$-term deviations. All commonly used dictionaries turn out not to be universal. In particular, the least  rational deviations in Hardy space $H^2$ do not  form certain strictly monotone null sequences.
There are no universal dictionaries in finite dimensional Hilbert spaces. We
  construct a universal dictionary in every  infinite dimensional Hilbert space.

\end{abstract}

\maketitle

\section{Introduction}

Let $H$ denote a real or complex Hilbert space with the norm $|\cdot|$ and the scalar product $\langle\cdot,\cdot\rangle$.
Let $D$ be a {\it dictionary}, that is, a
subset of the unit sphere $S(H)$ so that $\overline{\spn}D=H$.  For $m\in \NN$, $x\in H$ and
$$
\Sigma_m(D):=\left\{\sum_{k=1}^m \lambda_k g_k:\,  \lambda_k \in \R (\C),\,  g_k\in D\right\},
$$
$m$-term deviations are defined as
$$
\sigma_m(x)=\dist(x,\Sigma_m(D)).
$$
It is quite natural to define $\sigma_0(x)=|x|$.

In this general setting, $m$-term approximation has been introduced and first studied almost simultaneously by several authors~\cite{DT,DMA,DGDS}. Particular cases of approximation which occur to be $m$-term, such as  approximation by rational functions, splines, ridge functions, were widely investigated during the last century. The first of these cases seems to be that of bilinear approximation studied by E. Schmidt~\cite{S}; see Example~\ref{Schmidt}. S.B. Stechkin~\cite{St} was the first to  consider the setting of the problem of  $m$-term approximation with respect to a particular dictionary, the orthonormal basis in $H$.

In this paper, we study the set
$$
\sigma(D)=\{\{\sigma_m(x)\}_{m=0}^\infty:\,  x\in H\}
$$
of all sequences of $m$-term deviations for different dictionaries $D$.

In Theorem~\ref{theorem1} we prove a dichotomy result for $m$-term deviations: for any dictionary $D$, either
\begin{enumerate}
\item[(i)] $\sigma_m(x)$ decrease exponentially for every $x\in H$, or
\item[(ii)] for any sequence $\alpha_n\to 0$ there exists an element $x\in H$ such that $\sigma_m(x)\ge \alpha_m$, $m=0,1,2,\dots$
\end{enumerate}

A dichotomy in the rate of convergence  occurs in several approximation processes. For alternating projections onto a finite family of subspaces of a Hilbert space, it was  obtained independently in \cite{BDH,DH} and \cite{BaGM1,BaGM2}. In \cite{K} there is a parallel to this result concerning random products.  In~\cite{BK}, a dichotomy in the rate of convergence  was proved for greedy approximation with respect to an arbitrary dictionary in a Hilbert space. The property (ii) of arbitrarily slow convergence was established for different approximation schemes, in particular, for $m$-term approximation in Banach spaces with respect to dictionaries satisfying specific conditions~\cite{AO}.

Theorem 1 raises a natural question. Assume  the second  case (ii). For any strictly decreasing sequence $\alpha_0>\alpha_1>\dots>\alpha_m\to 0$, does there exist an element $x\in H$ with exact equalities $\sigma_m(x)=\alpha_m$, $m=0,1,2,\dots$? Note that we cannot expect to have such an element for a non-strictly decreasing $\alpha_m$: proximality of the set $\Sigma_m(D)$ in $H$ implies the strict inequality $\sigma_{m+1}(x)<\sigma_m(x)$ for any $x$.

This question is an analogue of the  {\em Bernstein  lethargy problem} for linear approximation: given a nested system $Y_1\subsetneq Y_2 \subsetneq \dots$ of linear subspaces of a Banach space $X$ and a strictly decreasing sequence $\alpha_0>\alpha_1>\dots>\alpha_m\to 0$, does there exist an element $x\in X$ with exact equalities $\dist (x,Y_m)=\alpha_m$, $m=0,1,2,\dots$? The answer is {\em yes} in many particular cases: if $X$ is a Hilbert space; if all $Y_m$ are finite-dimensional; if $\alpha_m>\sum_{n=m+1}^\infty \alpha_n$ for all $m$. However, the Bernstein problem is still unsolved in its general setting. Surveys of results related to this problem can be found in~\cite{B},~\cite{Kon}.

The question on the existence of an element with prescribed $m$-term deviations  can be treated as a nonlinear finite-dimensional version of the Bernstein problem. Surprisingly, we at once obtain  the negative  answer to this question for   commonly used  dictionaries satisfying the condition (ii). In particular, in Theorem \ref{theorem2} we prove that the least rational deviations cannot form an arbitrary strictly monotone sequence in the Hardy space $H^2$ in the upper half-plane.

So we have to modify the Bernstein-type problem mentioned above:
does there {\em at all exist} a dictionary $D\subset S(H)$ such that $\sigma(D)$ contains the entire set of strictly decreasing sequences tending to zero?  We call such a dictionary {\em universal}.

Finding a universal dictionary is challenging   for a  Hilbert space in particular. In the general Banach space setting, it is quite easy to give examples of universal dictionaries. For instance, the dictionary of the standard basis elements in $c_0$ is universal.

We construct a simple example of a universal dictionary in any non-separable Hilbert space in Remark~\ref{remark5}.  The rather complicated
example of a universal dictionary in a separable infinite dimensional Hilbert  space  we present in Theorem~\ref{theorem3}.

In Theorem~\ref{theorem4} we provide dictionaries $D$ in finite-dimensional Euclidean space, for which $\sigma(D)$ contains all not too slowly decreasing finite sequences.

\section{Dichotomy}

For a dictionary $D\subset S(H)$, we define
$$
\rho(D)=\inf_{x\in S(H)} \sup\{|\langle x,g\rangle| :\, g \in  D\},
$$
and a new increasing family of dictionaries
$D_m=\Sigma_m(D)\cap S(H)$, $m\in \NN$.
%Note that $D=D_1\subset D_2\subset D_3\subset \dots$.
The characteristic $\rho(D)$ influences the rate of convergence of the greedy algorithm.
 If $\rho(D)>0$, the algorithm converges fast everywhere; if $\rho(D)=0$,
 it converges arbitrarily slowly for certain starting elements \cite{BK}.
 In Theorem~\ref{theorem1} we will show a parallel to this result for the $m$-term approximation.
 Here the decisive property is whether $\rho(D_k)=0$ for all $k\in \NN$.

 Geometrically the condition   $\rho(D)=0$ means that  the dictionary $D$ is very ``slim": it is contained in an ``arbitrarily thin board" \cite{BK}. In particular, if  $\rho(D)=0$ then the interior of $\overline{\co}D$ is empty.

\begin{remark}\label{remark0}
Let $D\subset S(H)$ be a dictionary.
If $\overline{\co}D$ contains a ball of radius $r>0$,
then  $\rho(D)\geq r$.
\end{remark}
\begin{proof}
Assume $B(v,r)\subset \overline{\co}D$ and $x\in  S(H)$. Then
\begin{equation}\notag
\begin{split}
sup&\{|\langle x,g\rangle| :\, g \in  D\}=\sup\{|\langle x,g\rangle| :\, g \in  \overline{\co}D\} \\&\geq
\max\{|\langle x,v+rx\rangle|, |\langle x,v-rx\rangle|\}=\max\{|r\pm \langle x,v\rangle|\}\geq r.
\end{split}
\end{equation}

\end{proof}

\begin{lemmaX}\label{lemmaA} $($\cite{BK}$)$
Let $D\subset S(H)$ be a dictionary.
The equality $\rho(D)=0$ holds if and only if there exists an orthonormal sequence $\{w_n\}$ in $H$ so that
$$
\lim_{n\to \infty}\sup\{|\langle w_n,g\rangle| :\, g \in  D\}=0.
$$
\end{lemmaX}

\begin{theorem}\label{theorem1}
Let $D\subset S(H)$ be a dictionary.
\begin{enumerate}
\item[(i)]   If $\rho(D_k)>0$ for some $k\in \NN$, then
\begin{equation}\label{m-term-rate}
\sigma_m(x)\le |x| (1-\rho(D_k)^2)^{[m/k]/2},
\qquad m=0,1,2,\dots
\end{equation}
for every $x\in H$.
\item[(ii)]   If $\rho(D_k)=0$ for all $k\in \NN$, then  for every sequence $\alpha_m\to 0$ there exists  $x\in H$ such that  $\sigma_m(x)\geq \alpha_m$ for all $m=0,1,2,\dots $
\end{enumerate}
\end{theorem}
\begin{proof} (i) According to the definition  of $\rho(D_k)$, for any $x\in H$ we get
\begin{equation}\notag
\sigma_k^2(x)=|x|^2-\sup_{g\in D_k}|\langle x, g\rangle |^2\leq |x|^2(1-\rho(D_k)^2),
\end{equation}
so that for every $\varepsilon>0$ we can reduce the norm of $x$ with a coefficient of $(1-\rho(D_k)^2+\varepsilon)^{1/2}$ by subtracting an element from $\Sigma_k(D)$. In $m$-term approximation of $x$, we can subtract  $[m/k]$-times   and thus reduce the norm of $x$ with a coefficient of $(1-\rho(D_k)^2+\varepsilon)^{[m/k]/2}$. Since $\varepsilon$ is arbitrary, we get (\ref{m-term-rate}).

(ii) Since every null-sequence  admits a dominating decreasing null-sequence   we can assume $\alpha_0>\alpha_1>\dots$.
We can also assume that $\alpha_0\leq 1/12$:
if $\alpha_0>1/12$ and $x\in H$ works for the sequence $\{\alpha_m/(12\alpha_0)\}$, then $12\alpha_0x$ works for $\{\alpha_m\}$.
We define
$$
\beta_m=4\sqrt{\alpha_m^2-\alpha_{m+1}^2}.
$$
Then
$$
\sum_{n=m}^{\infty}\beta_n^2=16\alpha_m^2,
$$
hence $(\beta_m)\in \ell_2$.
In the spirit of Lemma~\ref{lemmaA}, we take an orthonormal  sequence $\{w_n\}$, such that for all $n\in \NN$
\begin{equation}\label{oncl}
\sup_{g\in D_n}|\langle w_n,g\rangle|\leq \beta_n/2.
\end{equation}
Here is an explanation, why does the orthonormal  sequence exist.
We choose a   sequence $v_k\in S(H)$ so that $\sup\{|\langle v_k,g\rangle| :\, g \in  D_k\}\leq \beta_k/4$ for all $k\in \NN$.  The sequence $\{v_k\}$ possesses a weakly convergent subsequence. This subsequence  converges weakly to zero, since $\overline{\spn} D=H$ and $D\subset D_k$ for all $k\in \NN$.   Hence there is a   an orthonormal sequence $\{w_n\}$ and a subsequence  $\{v_{k_n}\}$ of the above subsequence, so that
$|w_n-v_{k_n}|\leq \beta_n/4$ for all $n\in \NN$ (see {\em e.g.} Lemma~6.2 of~\cite{K}).   Since $D\subset D_2\subset D_3\subset\dots$, the estimates (\ref{oncl}) follow.

We define
$$
x=\sum_{n=0}^\infty\beta_nw_n.
$$
Then $\sigma_0(x)=|x|=4\alpha_0\leq 1/3$ and $|\sum_{n=m}^{\infty}\beta_nw_n|=4\alpha_m$.
Assume
$y\in \Sigma_m(D)$
is such that $|x-y|\leq 2\sigma_m(x)$. Then $|y|\leq 3|x|\leq 1$.
Hence for all $m\in \NN$ we have
\begin{equation}\notag
\begin{split}
\sigma_m(x)&\geq|x-y|/2\geq\frac 12\left|\left\langle x-y,\frac1{4\alpha_m} \sum_{n=m}^{\infty}\beta_nw_n\right\rangle\right|\\
&\geq\frac 12\left(\frac1{4\alpha_m} \sum_{n=m}^{\infty}\beta_n^2- \frac1{4\alpha_m} \sum_{n=m}^{\infty}\beta_n|\langle y,w_n\rangle | \right)\\
&\geq \frac 12\left(4\alpha_m-\frac1{4\alpha_m}\sum_{n=m}^{\infty}\beta_n^2/2\right)\\
&=\frac 12\left(4\alpha_m-\frac{16\alpha_m^2}{8\alpha_m}\right)=\alpha_m.
\end{split}
\end{equation}
\end{proof}

It is essential to require  $\rho(D_k)=0$ in Theorem~\ref{theorem1} for {\em all} $k\in \NN$ and not just for $k=1$.
Indeed, in Remark~\ref{remark1} we build for every  $m\in \NN$   a slim dictionary which stays slim being ``added" $(m-1)$-times to itself, but after adding the $m$-th copy it bloats.

\begin{remark}\label{remark1}
Let $H$ be an infinite dimensional Hilbert space. For every $m\in \NN$, $m\geq 2$ there exists a dictionary $D\subset S(H)$ so that $\rho(D)=\rho(D_2)=\dots=\rho(D_{m-1})=0$ and at the same time  $\rho(D_m)>0$.
\end{remark}
\begin{proof} Assume that $H$ is separable.
We write $H=\R^m_1\oplus\R^m_2\oplus\dots$ as a sum of $m$-dimensional orthogonal subspaces $\R^m_n$.
In each $\R^m_n$ we choose a unit vector $s_n$ and a basis of $m$ unit elements
$ g_1^n,\dots,g_m^n$
with the following property:  any  $(m-1)$-term linear combination $w$ of $ g_1^n,\dots,g_m^n$ (i.e., $w=\sum_{j\not= i}\lambda_jg^n_j$ for some $i\in \{1,\dots,m\}$) satisfies the inequality
$$
\left|\langle w, s_n\rangle\right|\leq \varepsilon_n|w|,
$$
where  $\varepsilon_n\to 0$ is a fixed sequence of positive numbers. This basis exists according to the first part of Construction~\ref{build} in the Appendix.

Let $D$ be the set of all norm-one elements of $H$ of the form $\sum_{n=1}^{\infty}\alpha_n g_{j_n}^n$, where $j_n\in\{1,\dots, m\}$ and $\alpha_n\in \R$.
Since $\Sigma_m(D)$ is a dense subset of $H$, $\rho(D_m)=1$. At the same time, for any element $w\in \Sigma_{m-1}(D)$ we have $|\langle w, s_n\rangle|\le \varepsilon_n |w|$, so
$$
\rho(D)=\rho(D_2)=\dots=\rho(D_{m-1})=0
$$
by Lemma \ref{lemmaA}.

If $H$ is not separable, we write $H=X\oplus_\perp Y$, where $X$ is an infinite dimensional separable
Hilbert space. In $X$ we choose a dictionary $D$ as above
and define the dictionary of $H$ as $\DD=D\cup S(Y)$. Then
$$
\DD_k=\Sigma_k(\DD)\cap S(H)=D_k\cup S(Y)\cup((\Sigma_{k-1}
(D)+Y)\cap S(H)),
$$
hence $D_k\cup S(Y)\subset \DD_k\subset (\Sigma_k(D)+Y)\cap S(H)$.

Since $D_m$ is dense in $S(X)$, it follows
from the first inclusion  that $\overline{\co}\DD_m$ contains a ball of radius $\sqrt2/2$. Hence
$\rho(\DD_m)\geq \sqrt2/2$ according to Remark~\ref{remark0}.

Let $1\leq k<m$ be given.
Using Lemma \ref{lemmaA} we choose an orthonormal sequence $\{w_n\}$ in $X$ so that
$$
D_k\subset K:=\{v\in H:\, |\langle v,w_n\rangle|\leq 1/n, n\in \NN\}.
$$
Every element $v\in (\Sigma_k(D)+Y)\cap S(H)$ can be written as $v=\alpha g+y$, where $\alpha\in [0,1]$, $g\in D_k$, and $y\in Y$. Hence $\DD_k\subset(\Sigma_k(D)+Y)\cap S(H)\subset K$ and   $\rho(\DD_k)=0$ by Lemma \ref{lemmaA}.
\end{proof}

\section{Seeking a universal dictionary}
\complain{Search for a universal dictionary}

The set $\sigma(D)$ of a dictionary $D$ consists of decreasing null-sequences. Here we  look into the size and the structure of  $\sigma(D)$. According to Theorem~\ref{theorem1} there is a dichotomy:  the set  $\sigma(D)$ either contains only very fast converging sequences, or  $\sigma(D)$ contains sequences converging arbitrarily slowly. We pay attention to the second case: we wonder if and when    $\sigma(D)$ contains  every strictly decreasing  null-sequence. Such a
dictionary  we call  universal.

In the next two canonical examples of dictionaries the sequence of $m$-term deviations is square-convex for each element, that is $\sigma_{n-1}^2-\sigma_n^2\geq \sigma_n^2-\sigma_{n+1}^2$. Consequently,  these dictionaries are not universal.

\begin{example}\label{ONbasis}
Let $H$ be a separable Hilbert space with an orthonormal basis $\{e_n\}_{n=1}^\infty$ and let  $D=\{e_n: n=1,2,\dots\}$.
\end{example}
It is easy to see that for each $x=\sum x_n e_n\in H$ we have
\begin{equation}\notag%\label{ONB_dev}
\sigma_m(x)=\left(\sum_{n=m+1}^\infty |x_{k(n)}|^2\right)^{1/2}
\end{equation}
where $k:\NN\to \NN$
is a permutation making the Fourier coefficients of $x$ monotonically decreasing: $|x_{k(1)}|\geq |x_{k(2)}|\geq \dots$

\begin{example}\label{Schmidt}
Let $H=L_2([0,1]^2)$ and  $D=\{\varphi(t)\psi(s): \varphi, \psi \in L_2[0,1]\} $.
\end{example}
It is well known that for each function $f\in H$,
\begin{equation}\notag%\label{Schmidt_dev}
\sigma_m(f)=\left(\sum_{n=m+1}^\infty s_n(f)^2\right)^{1/2}
\end{equation}
where $\{s_n(f)\}$ is the monotonically ordered sequence of $s$-numbers, the Neumann-Schatten numbers, of the operator $x(t)\to \int_0^1f(t,u) x(u)\,du$ acting in $L_2[0,1]$ (see e.g. \cite[Ch.\,1]{teml_book}).
\vspace{2mm}

Below we show that neither the dictionary consisting of step functions, nor that consisiting of  linear fractional  functions is universal.
The reason behind it is that a universal dictionary  $D$  satisfies a much stronger condition than
 $\rho(D_m)=0$ for all $m\in \NN$.
 Namely, each  $D_{m}$   has to be contained in an arbitrarily thin board defined by an element from the next iterate $D_{m+1}$.

\begin{remark}\label{remark2}
If $D$ is a universal dictionary in $H$, then
\begin{equation}\label{thin_boards}
\inf_{h\in D_{m+1}}\sup_{g\in D_{m}}|\langle h,g\rangle|=0 \mbox{ for each } m\in \NN.
\end{equation}
\end{remark}
\begin{proof}
Let $m\in \NN$ be given.
Since $D$ is universal, for every $\varepsilon>0$ there exists an element $x\in H$ having $\sigma_0(x)=|x|=1$, $\sigma_{m}(x)=1-\varepsilon$ and $\sigma_{m+1}(x)=\varepsilon$. The first  two equalities imply that $|\langle x, g\rangle|< \sqrt{2\varepsilon}$ for each $g\in D_{m}$. The third equality provides an element $h\in D_{m+1}$ so that $|x-h|<2\varepsilon$ and hence $|\langle h, g\rangle|< \sqrt{2\varepsilon}+2\varepsilon$ for each $g\in D_{m}$.
\end{proof}

In the proof of (\ref{thin_boards})  we have used only a very weak version of the universality of $D$:
for every $m\in \NN$ there is $x\in H$ which is very poorly approximated by $D_m$ and almost realised in $D_{m+1}$.  It would be of interest to know whether the condition (\ref{thin_boards})  implies the universality of a dictionary in a Hilbert space.
Note that all non-universal dictionaries in this paper lack the property (\ref{thin_boards}).

\begin{example}\label{piecewise_constant}
In $H=L_2[0,1]$ consider a dictionary $D$ consisting  of normalized step-functions $ a \chi_{[0,s]}+b \chi_{[s,1]}$, where $a,b\in \R$ and  $s\in [0,1]$. Then
\begin{enumerate}
\item[(1)]   $\inf_{h\in D_2}\sup_{g\in D_1}|\langle h,g\rangle|=0$ and
\item[(2)]   $\inf_{h\in D_m}\sup_{g\in D_{m-1}}|\langle h,g\rangle|\geq 1/\sqrt{m+1}$ for every $m\geq 3$,
\end{enumerate}
hence $D$ is not universal by Remark~\ref{remark2}.
\end{example}

The set $\Sigma_m$ here consists of  piecewise constant functions with at most $m+1$ pieces (splines of degree 0 with non-fixed breakpoints). It is well known that the $m$-term deviations in this case can have the rate much slower than exponential (see, {\em e.g.}, \cite{DeVore}), so this dictionary is of the type (ii) in Theorem~\ref{theorem1}.

To show (1) we define $h_n=\sqrt n\cdot\chi_{[1/2, 1/2+1/n]}\in D_2$ for $n\geq 4$. Since
$|g(t)|\leq 2$  for every $g\in D_1=D$ and $t\in [1/2,3/4]$, the equality (1) follows.

To show (2) let $h\in D_m$ taking values $a_0,\dots, a_{m}$ on the consecutive  intervals $I_0,\dots,I_m$ of lengths $\lambda_0,\dots, \lambda_{m}$ be given. Then $\sum_{j=0}^m a_j^2\lambda_j=1$. Without loss of generality we can assume that $a_1^2\lambda_1\geq 1/(m+1)$.  Define $g=(1/\sqrt{\lambda_1})\cdot\chi_{I_1}\in D_2$. Then $|\langle h,g\rangle|=|a_1|/\sqrt{\lambda_1}\cdot\lambda_1\geq 1/\sqrt{m+1}$.

\begin{example}\label{rational}
Let $H$ be a complex  Hilbert function space. For the dictionary $D$ we take all linear fractional  functions $r(t)=(at+b)/(ct+d)$, where $a,b,c,d\in \C$ are so that $r\in S(H)$.
\end{example}
For $m\in \NN$, the closure of $\Sigma_m(D)$ is equal to the set of all rational functions of degree at most $m$ from $H$. Thus the $m$-term deviations $\sigma_m(f)$ coincide with the least rational deviations $R_m(f)$ of order $m$.  These have not necessarily exponential rate in many particular spaces~\cite[Ch. 7, 10]{LGM},~\cite{D}. In the next section, we prove that the linear fractional functions $D$ do not form a universal dictionary in the Hardy space of functions analytic in a half-plane. We have chosen  this particular function space  as the scalar products of rational functions in it are easy to calculate.

\section{Rational deviations do not form certain strictly decreasing null sequences}\label{nec_cond}

\begin{remark}\label{remark3}

Let $D$ be a dictionary in a complex space $H$ so that $\lambda g\in D$ for any $g\in D$ and $\lambda \in \C$, $|\lambda|=1$.    Then the equality
\begin{equation}\label{thin_boards21}
\inf_{h\in D_2}\sup_{g\in D}|\langle h,g\rangle|=0,
\end{equation}
is equivalent to
\begin{equation}\label{thin_boards21_spec}
\inf_{a,b\in D}\sup_{g\in D}\left|\left\langle \frac{a-b}{|a-b|},g\right\rangle\right|=0.
\end{equation}

\end{remark}

\begin{proof} Clearly, (\ref{thin_boards21_spec}) implies (\ref{thin_boards21}).

Let (\ref{thin_boards21}) hold, that is, for every $\varepsilon\in (0,1/2)$ there exists $h\in D_2$ so that $|\langle h,g\rangle|<\varepsilon$ for any $g\in D$. We may assume $h=(a-\lambda b)/\delta$, where $a,b\in D$, $\lambda\in [0,1]$, and $\delta=|a-\lambda b|$. We have
$$
\delta^2=|a-\lambda b|^2=1+\lambda^2-2\lambda \Re \langle a,b \rangle
$$
and
$$
1-\lambda \Re \langle a,b \rangle\le |1-\lambda \langle b,a \rangle|=|\langle a- \lambda b, a\rangle|\le \varepsilon\delta,
$$
so that
$$
1-\frac{1+\lambda^2-\delta^2}{2}\le \varepsilon\delta \Longrightarrow 2\varepsilon\delta\ge 1-\lambda^2+\delta^2> 1-\lambda.
$$
Consequently,
$$
\left|\left\langle  \frac{a-b}{|a-b|},g \right\rangle\right|\le \frac{|\langle a- \lambda b, g\rangle| +(1-\lambda)}{|a- \lambda b|-(1-\lambda)} \le \frac{\varepsilon\delta+2\varepsilon\delta}{\delta-2\varepsilon\delta}=\frac{3\varepsilon}{1-2\varepsilon},
$$
so that (\ref{thin_boards21_spec}) holds.
\end{proof}

We consider the Hardy space $H^2(\Pi_+)$ of complex functions $f(z)$ holomorphic in the upper half-plane $\Pi_+=\{z: \Im z>0\}$, for which
$$
\|f\|:=\sup_{y>0} \left(\int_{\R} |f(x+iy)|^2\, dx \right)^{1/2} < \infty.
$$
It is well known that each function $f\in H^2(\Pi_+)$ has angular limits $f(x)$ for almost all $x\in \R$, and the $L_2(\R)$-norm of this limit function coincides with $\|f\|$ (see {\em  e.g.}~\cite[Ch.\,6]{Koo}).

In $H^2(\Pi_+)$ the dictionary from Example \ref{rational} consists of normalized functions
  of the form
$$
\frac{c}{z-a}, \qquad a\in \Pi_-= \{z: \Im z<0\}, c\in \C.
$$

Since
$$
\left\|\frac{1}{z-a}\right\|^2=\int_{\R}\frac{dx}{(x-a)(x-\bar{a})}=\frac{\pi}{|\Im a|},
$$
we are dealing with the dictionary
$$
D=\left\{\frac{e^{i\theta}}{\sqrt{\pi}}\frac{\sqrt{|\Im a|}}{z-a}: a\in \Pi_-, \theta\in \R\right\}.
$$
Clearly, $m$-term deviations with respect to  $D$ coincide with the least rational deviations in $H^2(\Pi_+)$.  This is the case even for $m=0$, as the only constant function in $H^2(\Pi_+)$ is identically zero.  The above dictionary $D$ is of the type (ii) in Theorem~\ref{theorem1} as
the rational deviations  can have their rate much slower than exponential \cite{Stel}.

\begin{theorem}\label{theorem2}
For this dictionary $D$,
\begin{equation}\label{not_thin_boards_rat}
\inf_{g_1,g_2\in D}\sup_{g\in D}\left|\left\langle \frac{g_1-g_2}{\|g_1-g_2\|},g\right\rangle\right|> 10^{-2}.
\end{equation}
By Remarks \ref{remark2} and \ref{remark3}, this means that $D$ is not universal in $H^2(\Pi_+)$.
\end{theorem}
\begin{proof}
1. We may assume
$$
g_1(z)=\frac{1}{\sqrt{\pi}} \frac{\sqrt{\alpha}}{z+\alpha i}, \qquad g_2(z)=\frac{1}{\sqrt{\pi}} \frac{\sqrt{\beta}e^{i\theta}}{z+\beta i -\gamma}
$$
where the poles $a=-\alpha i$, $b=-\beta i+ \gamma$ are such that $\alpha\ge \beta>0$, $\gamma\ge 0$, and $\theta\in [-\pi,\pi]$.

Next, we calculate the norm of $r=g_1-g_2$:
$$
\|r\|^2=\|g_1\|^2+\|g_2\|^2-2\Re\langle g_1, g_2 \rangle
$$
$$
=2-2\Re \frac{1}{\pi}\int_{\R}\frac{\sqrt{\alpha}\sqrt{\beta}e^{i\theta}}{(x-\bar{a})(x-b)}\, dx=2-\frac{2\sqrt{\alpha\beta}}{\pi}\Re\frac{2\pi ie^{i\theta}}{\bar{a}-b}.
$$
In the notation $\bar{a}-b=|\bar{a}-b|ie^{i\varphi}$, this is equal to
\begin{equation}\label{norm of r}
\begin{array}{l}
\|r\|^2=2-\sqrt{|\bar{a}-b|^2-|a-b|^2} \Re \frac{2e^{i(\theta-\varphi)}}{|\bar{a}-b|}\\
\qquad =2\left(1-\sqrt{1-\left|\frac{a-b}{\bar{a}-b}\right|^2}\cos(\theta-\varphi)\right).
\end{array}
\end{equation}

We aim to prove that
\begin{equation}\label{aim}
\sup_{g\in D}\left|\left\langle \frac{g_1-g_2}{\|r\|},g\right\rangle\right|> 10^{-2}.
\end{equation}

2. Clearly,
$$
\left\langle \frac{g_1-g_2}{\|r\|},g_1\right\rangle+\left\langle \frac{g_1-g_2}{\|r\|},-g_2\right\rangle=\|r\|.
$$

If $\|r\|>2\cdot 10^{-2}$, then the modulus of one of the terms is greater than $10^{-2}$, and (\ref{aim}) holds.

Hence we may assume $0<\|r\|\le 2\cdot 10^{-2}$. Together with (\ref{norm of r}) this implies
\begin{equation}\label{delta}
 1-\left|\frac{a-b}{\bar{a}-b}\right|^2\ge (1-2\cdot 10^{-4})^2 \Longrightarrow \left|\frac{a-b}{\bar{a}-b}\right|< 2\cdot 10^{-2},
\end{equation}
and
\begin{equation}\label{cos}
\cos(\theta-\varphi)> 1-2\cdot 10^{-4}.
\end{equation}

Recall that
$$
e^{i\varphi}=\frac{\bar{a}-b}{|\bar{a}-b|i} \Longrightarrow \sin \varphi =-\Re \frac{\bar{a}-b}{|\bar{a}-b|}=-\Re \frac{a-b}{|\bar{a}-b|},
$$
so that
\begin{equation}\label{sin}
0\le\sin \varphi\le \left|\frac{a-b}{\bar{a}-b}\right| < 2\cdot 10^{-2}.
\end{equation}
We use that $2x/\pi\leq \sin x$  for small  positive $x$ to get
 $\varphi\in [0,\pi 10^{-2}]$ and  $|\theta-\varphi|<\pi(1+10^{-2})$. This and (\ref{cos}) using that    $\cos x\leq 1-x^2/\pi$  for small   $x$ implies
\begin{equation}\label{theta}
1-|\theta-\varphi|^2/\pi> 1- 2\cdot 10^{-4} \Longrightarrow |\theta-\varphi|<3\cdot 10^{-2} \Longrightarrow |\theta|<7\cdot 10^{-2}.
\end{equation}

Now we can give an upper bound for the value of $\|r\|^2$ from (\ref{norm of r}). We use the notation $\left|\frac{a-b}{\bar{a}-b}\right|=\delta$, inequalities (\ref{sin}) and (\ref{theta}):
\begin{equation}\label{bound for norm of r}
\begin{array}{l}
\|r\|^2=2(1-\sqrt{1-\delta^2}(1-2\sin^2((\theta-\varphi)/2))\\
\qquad \le 2(1-(1-\delta^2)(1-2\sin^2((\theta-\varphi)/2))\\
\qquad \le 2(\delta^2+2\sin^2((\theta-\varphi)/2))\le 2(\delta^2+2(|\sin \frac{\theta}{2}| + |\sin \frac{\varphi}{2}|)^2)\\
\qquad \le 2(\delta^2+2(|\sin \theta| + |\sin \varphi|)^2)\le 2(\delta^2+4\sin^2 \theta + 4\sin^2 \varphi)\\
\qquad \le 2(5\delta^2+4\sin^2 \theta)=10\delta^2+8\sin^2 \theta.
\end{array}
\end{equation}

3. We consider the scalar product of $r/\|r\|$ and of an arbitrary element $g\in D$. Since $g$ has the form
$g(x)=e^{i\tau}\sqrt{|\Im w|}/(\sqrt{\pi}(x-w))$,
$$
\left\langle\frac{r}{\|r\|}, g \right\rangle
=\frac{e^{-i\tau} \sqrt{|\Im w|}}{\sqrt{\pi}\|r\|}\int_{\R}\frac{r(x)}{x-\bar{w}}=\frac{e^{-i\tau} \sqrt{|\Im w|}}{\sqrt{\pi}\|r\|}2\pi i r(\bar{w}),
$$
so that
$$
\left|\left\langle\frac{r}{\|r\|}, g\right\rangle\right|^2=\frac{4\pi |\Im z| |r(z)|^2}{\|r\|^2}=:f(z)
$$
for some $z=\bar{w}\in \Pi_+$.

Calculating
\begin{eqnarray*}
r(z)&=&\frac{1}{\sqrt{\pi}}\left(\frac{\sqrt{\alpha}}{z-a}-\frac{\sqrt{\beta}e^{i\theta}}{z-b}\right)\\
&=& \frac{(z-a)(\sqrt{\alpha}-\sqrt{\beta}e^{i\theta}) + \sqrt{\alpha}(a-b)}{\sqrt{\pi}(z-a)(z-b)} \\
&=&\frac{\sqrt{\alpha}((z-a)(1-\sqrt{\beta/\alpha}e^{i\theta}) + (a-b))}{\sqrt{\pi}(z-a)(z-b)},
\end{eqnarray*}
and using the bound (\ref{bound for norm of r}) for $\|r\|^2$, we get
\begin{equation}\label{f-bound}
f(z)\ge \frac{2\pi |\Im z| \alpha (|z-a||1-\sqrt{\beta/\alpha}e^{i\theta}|-|a-b|)^2}{|z-a|^2 |z-b|^2(5\delta^2+4\sin^2\theta)}.
\end{equation}
To prove (\ref{aim}) and consequently also (\ref{not_thin_boards_rat}), we have to find $z\in \Pi_+$, so that $f(z)\ge 10^{-4}$. The choice of $z$ depends on the relations between $\theta$, $\delta$, $a$ and $b$.

4. In the case when $|\sin \theta| \ge \delta$, we have
$$
|a-b|\le |\bar{a}-b| |\sin \theta|\le ( |\bar{a}-a|+|a-b|) |\sin \theta|,
$$
and then using  (\ref{theta})
$$
|a-b|\le \frac{|\bar{a}-a| |\sin \theta|}{1-|\sin \theta|}\le \frac{2\alpha |\sin \theta|}{1-7\cdot 10^{-2}}< 3\alpha |\sin \theta|.
$$
This implies
$$
|3\alpha i-b|=|3\bar{a}-b|\le |3\bar{a}-a|+|a-b|\le 4\alpha+3\alpha|\sin \theta|< 5\alpha,
$$
$$
\Longrightarrow |(3\alpha i-a)(1-\sqrt{\beta/\alpha}e^{i\theta})|=4\alpha|1-\sqrt{\beta/\alpha}e^{i\theta}|\ge 4\alpha|\sin \theta|,
$$
$$
5\delta^2+ 4\sin^2 \theta\le 9 \sin^2 \theta,
$$
so that, according to (\ref{f-bound}),
$$
f(3\alpha i)\ge \frac{2\pi 3\alpha^2 (|(3\alpha i-a)(1-\sqrt{\beta/\alpha}e^{i\theta})|-|a-b|)^2}{|3\alpha i-a|^2 |3\alpha i-b|^2(5\delta^2+4\sin^2\theta)}
$$
$$
\ge \frac{6\pi \alpha^2 (4\alpha|\sin \theta|-3\alpha|\sin \theta|)^2}{(4\alpha)^2(5\alpha)^2 9\sin^2 \theta }=\frac{\pi}{600}>10^{-4}.
$$

5. In the case when $|\sin \theta| < \delta$, we have
$$
5\delta^2+ 4\sin^2 \theta< 9 \delta^2.
$$

If, in addition, $|1-\sqrt{\beta/\alpha}e^{i\theta}|\le\delta/3$, then
$$
|(\alpha i-a)(1-\sqrt{\beta/\alpha}e^{i\theta})|=2\alpha|1-\sqrt{\beta/\alpha}e^{i\theta}|\le 2|\bar{a}-b|\frac{\delta}{3}=\frac{2}{3}|a-b|,
$$
so that, according to (\ref{f-bound}),
$$
f(\alpha i=\bar{a})\ge \frac{2\pi \alpha^2 (|(\alpha i-a)(1-\sqrt{\beta/\alpha}e^{i\theta})|-|a-b|)^2}{(2\alpha)^2 |\bar{a}-b|^29\delta^2}
$$
$$
\ge \frac{\pi(|a-b|/3)^2}{18|a-b|^2}=\frac{\pi}{162}>10^{-4}.
$$

If, alternatively, $|1-\sqrt{\beta/\alpha}e^{i\theta}|>\delta/3$, then
$$
|(8\alpha i-a)(1-\sqrt{\beta/\alpha}e^{i\theta})|=9\alpha|1-\sqrt{\beta/\alpha}e^{i\theta}|>3\alpha\delta=\frac{3}{2} |\bar{a}-a|\delta
$$
$$
\ge \frac{3}{2} (|\bar{a}-b|-|a-b|)\delta= \frac{3}{2} |a-b|(1-\delta)>\frac{4}{3}|a-b|
$$
-- in view of (\ref{delta}),
$$
|8\alpha i-b|=|8\bar{a}-b|\le |\bar{a}-b|+\frac{7}{2}|\bar{a}-a|\le |\bar{a}-b|+\frac{7}{2}|\bar{a}-b|+\frac{7}{2}|a-b|
$$
$$
=\frac{9}{2}|\bar{a}-b|+\frac{7}{2}|\bar{a}-b|\delta< 5|\bar{a}-b|,
$$
so that, according to (\ref{f-bound}),
$$
f(8\alpha i)\ge \frac{2\pi 8\alpha^2 (|(8\alpha i-a)(1-\sqrt{\beta/\alpha}e^{i\theta})|-|a-b|)^2}{(9\alpha)^2 |8\alpha i-b|^29\delta^2}
$$
$$
\ge \frac{16\pi(|a-b|/3)^2}{81\cdot 25|\bar{a}-b|^29\delta^2}=\frac{16\pi}{164025}>10^{-4}.
$$
\end{proof}

Apparently, it is proved here for the first time that the least  rational deviations cannot form an arbitrary strictly monotone sequence in the Euclidean norm. In this sense, the uniform norm is better: A.A. Pekarskiĭ~\cite{Pek} proved that any strictly monotone sequence realizes as the sequence of the least rational deviations in the space  $C^{\C}[0,1]$ of complex continuous functions with the uniform norm.

In Theorem~\ref{theorem2} we have shown that monotone sequences with large jumps at the beginning  cannot be realized as sequences of rational deviations.
We leave open the question if there is  $k\in \NN$ so that  for any strictly decreasing sequence $\{\alpha_m\}$ there  exists a function $f\in H^2$ with   $\sigma_m(f)=R_m(f)=\alpha_m$ for all $m>k$.

\section{Non-separable case}

\begin{remark}\label{remark4}
For every sequence $d_n \downarrow 0$ of positive numbers, there exists a dictionary $D\subset l_2$ and an element $x\in l_2$ such that
 $\sigma_n(x)=d_n$, $n=0,1,2,\dots$.
\end{remark}
\begin{proof}
Let $e_0, e_1,\dots$ be  a standard basis in $l_2$. We set
$$
x=\sum_{n=0}^\infty\sqrt{d_n^2-d_{n+1}^2}e_n
$$
and
$$
\begin{array}{l}
g_0=e_0,\\
g_1=\sqrt{d_1^2-d_2^2}e_0-\sqrt{d_0^2-d_1^2}e_1,\\
g_2=\sqrt{d_2^2-d_3^2}e_1-\sqrt{d_1^2-d_2^2}e_2,\\
........\\
g_n=\sqrt{d_n^2-d_{n+1}^2}e_{n-1}-\sqrt{d_{n-1}^2-d_n^2}e_n,\\
........
\end{array}
$$
Since $d_n$ is strictly decreasing to zero,
$$
\spn\{g_0,g_1,\dots, g_n\}=\spn\{e_0,e_1,\dots, e_n\}$$
and $\overline{\spn}\{g_0,g_1,\dots\}=H$.  Hence
  $D=\{g_n/|g_n|\}_0^\infty$ is a dictionary. Note that $x$ is orthogonal to each of $g_1,g_2,\dots$

Now we calculate $n$-term deviations $\sigma_n(x)$ with respect to $D$. Clearly, $\sigma_0(x)=|x|=d_0$. Next,
$$
\dist(x,\spn\{g_0,\dots,g_{n-1}\})=\dist(x,\spn\{e_0,\dots,e_{n-1}\})=d_n,
$$
so that $\sigma_n(x)\leq d_n$, $n=1,2,\dots$.

Given any $n$-tuple $0\leq i_1<i_2<\dots<i_n$, let $k$ be the smallest non-negative
 integer which does not belong to $\{i_1,\dots,i_n\}$.  Three cases are possible.

If $k=0$, then $x\perp \spn\{g_{i_1},\dots,g_{i_n}\}$, and
$$
\dist(x,\spn\{g_{i_1},\dots,g_{i_n}\})=|x|=d_0>d_n.
$$

If $k=n$, then $\{i_1,\dots,i_n\}=\{0,\dots,n-1\}$, and
$$
\dist(x,\spn\{g_{i_1},\dots,g_{i_n}\})=d_n.
$$

If $i_j<k<i_{j+1}$ for some $j$, then $i_1=0,\dots,i_j=j-1$, $i_{j+1}\geq j+1$, so that
$$
\spn\{g_{i_1},\dots,g_{i_n}\} = \spn\{g_0,\dots,g_{j-1}\} \oplus_{\perp} \spn\{g_{i_{j+1}},\dots,g_{i_n}\},
$$
and $x\perp \spn\{g_{i_{j+1}},\dots,g_{i_n}\}$. Consequently,
$$
\dist(x,\spn\{g_{i_1},\dots,g_{i_n}\})=\dist(x,\spn\{g_0,\dots,g_{j-1}\})=d_j>d_n.
$$

We conclude that $\sigma_n(x)= d_n$.

\end{proof}

Remark \ref{remark4} can be used to construct a universal dictionary in a non-separable Hilbert space. This construction may seem superfluous in light of the fact that we intend to build a universal dictionary in a separable Hilbert space in the next section, and this dictionary can be easily extended to a universal dictionary in the non-separable case. Indeed, suppose $D$ is a universal dictionary of a Hilbert space $X$. If $H=X\oplus_\perp Y$,   then  $\DD=D\cup S(Y)$  is a universal dictionary in $H$: for every element in $X$ the best $m$-term approximation stays in $X$. However,  we present  the following direct construction because it is much simpler than that of Theorem \ref{theorem3}.

\begin{remark}\label{remark5}

In the non-separable Hilbert space $l_2([0,1]\times \{0,1,2,\dots\})$, it is quite easy to construct a universal dictionary.

\end{remark}

\begin{proof}
Let $\alpha:t\mapsto \alpha(t)$ be a bijection between $[0,1]$ and the set of all positive sequences in the usual space $l_2$, $\alpha(t)=(\alpha^t_0,\alpha^t_1,\dots)$ ($\alpha^t_n>0$ for all $t$ and $n$).
Let $\{e^t_n\}_{n=0}^\infty$ be an orthonormal basis in $l_2(\{t\}\times \{0,1,2,\dots\})$, $t\in [0,1]$.
We set
$$
\begin{array}{l}
g^t_0=e^t_0,\\
g^t_1=\alpha^t_1e^t_0-\alpha^t_0e^t_1,\\
g^t_2=\alpha^t_2e^t_1-\alpha^t_1e^t_2,\\
........\\
g^t_n=\alpha^t_ne^t_{n-1}-\alpha^t_{n-1}e^t_n,\\
........
\end{array}
$$
It is clear that $D=\{g^t_n/|g^t_n|: t\in [0,1], n=0,1,\dots\}$ is a dictionary. By the proof of Remark \ref{remark4}, the element
$$
x^t=\sum_{n=0}^\infty\alpha^t_ne^t_n
$$
has $\sigma_m(x^t)=\sum_{n=m}^\infty(\alpha^t_n)^2$ ($m=0,1,\dots$) for each $t\in [0,1]$, so that every strictly monotonic sequence of deviations is realized.
\end{proof}

\section{Universal dictionary in a separable Hilbert space}

\begin{theorem}\label{theorem3}
In a separable infinite dimensional Hilbert space $H$, there exists a dictionary $D$ such that for any sequence $d_n\to 0$, either strictly monotonic ($d_0>d_1>d_2>\dots$) or strictly monotonic down to zero ($d_0>d_1>\dots>d_N=0=d_{N+1}=d_{N+2}=\dots$)
there exists an element $x\in H$ having $\sigma_m(x)=d_m$, $m=0,1,2,\dots$
\end{theorem}
\begin{proof}

1. In the sequence space $l_2$ we consider an orthonormal basis enumerated as $\{e_0,e_{n,k}: n\in \NN, k\in \NN\}$. Let $\SSS=\{s_\nu\}_{\nu=1}^\infty\subset l_2$ be the countable dense set of all finitely supported non-zero elements with rational coordinates:
each $s_\nu$ has  the coordinates
$$
s_\nu=(s^\nu_0,s^\nu_{1,1}, s^\nu_{1,2},\dots, s^\nu_{n,k}, \dots ).
$$
We denote
\begin{equation}\notag
\begin{split}
n(\nu)&=\max\{n:\,  s^\nu_{n,k}\not=0 \  \mbox{for some}\ k\} \\
k(\nu)&=\max\{k:\,  s^\nu_{n,k}\not=0\  \mbox{for some}\ n\}\\
\eps(\nu)&=1/k(\nu)
\end{split}
\end{equation}
If $s_\nu\in \spn\{e_0\}$, we define $n(\nu)=k(\nu)=0$.

2. For each $\nu\in \NN$, let     $H_\nu$  be a $n(\nu)$-dimensional Euclidean space, and let
$H$ be the separable Hilbert space
$$
H=l_2\oplus_\perp H_1\oplus_\perp H_2\oplus_\perp\dots.
$$
According to Construction~\ref{build} in the Appendix,
in each subspace  $\spn\{s_\nu\}\oplus_\perp H_\nu$ of $H$ with $n(\nu)>0$ there is  a basis of $n(\nu)+1$ unit vectors $g_0^\nu,\dots,g_{n(\nu)}^\nu$
with the following property. Any  $m$-term linear combination $w$ of $g_j^\nu$'s with all the relevant
$n(\nu)\geq m$  and all the relevant $\eps(\nu)<\eps$ satisfies the inequality
\begin{equation}\label{proj}
|Pw|\leq \varepsilon\sqrt m|(Id-P)w|;
\end{equation}
here $P:H\to \ell_2$ is the orthogonal projection.

The span of the set
$$
D=\{e_0, g_i^\nu:\ \nu\in \NN,\ i=0,1,\dots,n(\nu)\}.
$$
contains the spanning set $\{e_0, s_\nu:\, \nu\in \NN\}$ of $\ell_2$ and   each $H_{\nu}$, hence  $D$ is a dictionary in $H$.

3. Let $k_1<k_2<\dots$ be any fixed increasing sequence of positive integers. We will prove that for each element of the form
$$
x= \alpha_0e_0+\sum_{n=1}^\infty \alpha_n e_{n,k_n},
$$
the $m$-term deviations with respect to $D$ can be estimated as
\begin{equation}\label{1dev_estim}
\sqrt{1-\frac{m}{k_m^2}} \left(\sum_{n=m}^\infty \alpha_n^2\right)^{1/2}\leq \sigma_m(x)\leq \left(\sum_{n=m}^\infty \alpha_n^2\right)^{1/2}
\end{equation}
for each $m=0,1,2,\dots$

4. We denote
\begin{equation}\notag
\begin{split}
E_m&=\spn\{e_0, e_{1,k_1},\dots, e_{m-1,k_{m-1}}\}, \\
F_m&=\spn\{e_0, e_{n,k}:\, n<m \mbox{ or } k < k_m\}.
\end{split}
\end{equation}
By $P_m$ we denote the orthogonal projection of $H$ onto $F_m$ and $x_m=P_mx\in E_m$, $r_m=(Id-P_m)x$.
Since the sequence $\{k_n\}$ is increasing,
$$
x_m= \alpha_0e_0+\sum_{n=1}^{m-1} \alpha_n e_{n,k_n},\ r_m=\sum_{n=m}^\infty \alpha_n e_{n,k_n},\, |r_m|^2=\sum_{n=m}^\infty \alpha_n^2
$$
Given any $\varepsilon>0$, there exists  $s_\nu\in E_m$ so that
$|x_m-s_\nu|<\varepsilon$. Hence $n(\nu)\leq m-1$, and
$$
\sigma_m(x)\leq \dist (x, \spn\{g_0^\nu,\dots,g_{n(\nu)}^\nu\})\leq |x-s_\nu|<\left(\sum_{n=m}^\infty \alpha_n^2\right)^{1/2}+\varepsilon,
$$
and the upper bound in (\ref{1dev_estim}) follows.

5. Let $y\in \Sigma_m(D)$, $y=y_1+\dots+y_l$, where $l\leq m$ and $y_j\in \spn\{s_{\nu_j}\}\oplus_\perp H_{\nu_j}$ are linear combinations of not more than $q_j$ elements of $D$, $q_1+\dots+q_l\leq m$.  We divide the summands into two groups: those with $s_{\nu_j}\in F_m$ and the rest. Adding up the summands within each of the two groups we get
$y=z+w$ with $(Id-P_m)Pz=0$ in the first group. In the second group all $n(\nu_j)\geq m$ and   $k(\nu_j)\geq k_m$.  This means in turn that in the second group $\eps(\nu_j)\leq 1/k_m$, hence according to (\ref{proj})
\begin{equation}\label{projw}
|Pw|\leq  \sqrt m|(Id-P)w|/k_m= \sqrt m R/k_m;
\end{equation}
Here we have denoted $R=|(Id-P)w|$
for brevity. Using the Pythagorean theorem several times and (\ref{projw}) toward the end we estimate
\begin{equation}\notag
\begin{split}
|x-y|^2&=|x-Py|^2+|(Id-P)y|^2 \\
&=|x-Pz-Pw|^2+|(Id-P)z|^2+|(Id-P)w|^2 \\
&\geq |x-Pz-Pw|^2+ R^2 \\
&=|P_m(x-Pz-Pw)|^2+|(Id-P_m)(x-Pz-Pw)|^2+ R^2 \\
&\geq|(Id-P_m)(x-Pz-Pw)|^2+ R^2 \\
&=|r_m-(Id-P_m)Pw|^2+ R^2 \geq  |r_m|^2-2|r_m||(Id-P_m)Pw|  + R^2 \\
&\geq  |r_m|^2-2|r_m||Pw|  +R^2\geq  |r_m|^2-2 \sqrt m|r_m|R/k_m+ R^2 \\
&=|r_m|^2(1-m/k_m^2)+(|r_m|\sqrt m/k_m-R)^2 \\
&\geq|r_m|^2(1-m/k_m^2).
\end{split}
\end{equation}
This proves the lower bound of (\ref{1dev_estim}), as $|r_m|^2=\sum_{n=m}^\infty \alpha_n^2$.

6. Let now a strictly decreasing sequence $d_n\downarrow 0$ be given. We set
\begin{equation}\label{eps}
\varepsilon_n=\min\{d_0-d_1,d_1-d_2,\dots, d_{n-1}-d_n\}, \qquad n=1,2,\dots
\end{equation}
We choose a strictly increasing sequence of positive integers $k_n$ so that
\begin{equation}\label{k_n}
k_n\geq \sqrt{\frac{d_0n}{\varepsilon_n}}, \qquad n=1,2,\dots,
\end{equation}
and consider the mapping $t\mapsto x(t)$ of the compact set
$$
K_\varepsilon=\{t=(t_1,t_2,\dots )\in c_0: 0\leq t_n\leq \varepsilon_n\}\subset c_0
$$
into $H$:
$$
 x(t)=\sqrt{d_0^2-(d_1+t_1)^2}e_0+\sum_{n=1}^\infty \sqrt{(d_n+t_n)^2-(d_{n+1}+t_{n+1})^2}e_{n,k_n}.
$$
This mapping is, due to (\ref{eps}), well defined and continuous.   From  (\ref{1dev_estim}) it follows  that for all $m\in \NN$
$$
\sigma_m(x(t))\leq d_m+t_m,
$$
and using also (\ref{k_n})
$$
\sigma_m(x(t))\geq \sqrt{1-\frac{m}{k_m^2}} (d_m+t_m)\geq \left(1-\frac{m}{k_m^2}\right) (d_m+t_m)
$$
$$
\geq d_m+t_m-\frac{d_0m}{k_m^2}\geq d_m+t_m-\varepsilon_m.
$$
Thus, the mappings $t\mapsto x(t)$ and $\sigma: x\mapsto\{\sigma_m(x)\}$ satisfy all the assumptions of

\begin{lemmaX}\label{lemmaB} $($\cite{Pek}$)$
 Let $\sigma:H\to c_0$, $x\to (\sigma_1(x), \sigma_2(x), \dots)\in c_0$ be a continuous mapping.
Let $\{d_n\}_{n=1}^\infty$ be a sequence of positive numbers  strictly decreasing  to zero.
Let $\{\varepsilon_n\}_{n=1}^\infty$ be a  sequence of non-negative numbers tending to zero and let
$$
K_\varepsilon=\{t=(t_1,t_2,\dots )\in c_0: 0\leq t_n\leq \varepsilon_n\}
$$
be a convex compact set in $c_0$ corresponding to this sequence.

If there is a continuous map $x:K_\varepsilon \to H$ such that
$$
d_n+t_n-\varepsilon_n\leq \sigma_n(x(t))\leq d_n+t_n, \qquad t\in K_\varepsilon, n=1,2,\dots,
$$
then there exists $x\in H$ such that $\sigma_n(x)=d_n$, $n=1,2,\dots$.
\end{lemmaX}

\begin{proof}
The mapping $F: K_\varepsilon\to K_\varepsilon$, $[F(t)]_n=d_n+t_n-\sigma_n(x(t))$, has a fixed point $t^*$ by Schauder theorem. It is clear that $x(t^*)$ is the element needed.
\end{proof}

In our case, Lemma \ref{lemmaB} provides an element $x=x(t^*)$ having $\sigma_m(x)=d_m$, $m=1,2,\dots$. Since $\sigma_0(x(t))=d_0$ for all $t$, we get the element required.

7. In case when a sequence strictly monotonic down to zero is given ($d_0>d_1>\dots>d_N=0=d_{N+1}=d_{N+2}=\dots$), we just repeat the arguments of part 6, having $\varepsilon_n=0$ for $n>N$ and the compact set $K_\varepsilon$ being finitely supported in $c_0$.

\end{proof}

\section{Finite-dimensional case}

In the Euclidean space $\R^N$, any dictionary satisfies condition (i) of Theorem \ref{theorem1} and hence is not universal. We present dictionaries realizing certain classes of $m$-term deviation sequences.

\begin{theorem}\label{theorem4}
For every $\varepsilon>0$, every $M>0$ and every $N\in \NN$, there exists a dictionary $D=D(\varepsilon,M,N)\subset \R^{N+1}$ with the following property. For every sequence
$M>d_0>d_1>\dots>d_N>0$ with $\min \{d_n-d_{n-1}: n=1,\dots,N\}>\varepsilon$ there exists an element $x\in \R^{N+1}$ having $\sigma_n(x)=d_n$, $n=0,\dots,N$.
\end{theorem}
\begin{proof}
Let $e_1,\dots,e_{N+1}$ be the standard basis of  $\R^{N+1}$.

According to Construction~\ref{matrix} in the Appendix there is another normalised basis $g_1,\dots,g_{N+1}$ with the following property.
For every $n\in \{1,\dots, N\}$, each $n$-dimensional subspace $\spn\{g_{i_1},\dots,g_{i_n}\}$ is so close to $\spn\{g_1,\dots,g_n\}=\spn\{e_1,\dots,e_n\}$ that the Hausdorff distance between their intersections with the ball $\{x\in \R^{N+1}:\, |x|\le M\}$ is less than $\varepsilon$. We define the desired  dictionary $D$ as the new basis:   $D=\{g_1,\dots,g_{N+1}\}$. The mapping
$t=(t_1,\dots,t_N)\to x(t)$ where
$$
x(t)=\sqrt{d_0^2-(d_1+t_1)^2} e_1+\sum_{k=2}^{N} \sqrt{(d_{k-1}+t_{k-1})^2-(d_k+t_k)^2} e_k+ (d_N+t_N)e_{N+1}
$$
maps the compact set $K_\varepsilon=\{t=(t_1,\dots t_N)\in \R^N: 0\leq t_n\leq \varepsilon\}$ continuously to   $\R^{N+1}$.  Moreover,
\begin{equation}\notag
\begin{split}
\sigma_n(x(t))&\leq \dist(x(t),\spn\{g_1,\dots,g_n\})=\dist(x(t),\spn\{e_1,\dots,e_n\}) \\
&=\left(\sum_{k=n+1}^{N+1}x(t)_k^2\right)^{1/2}=d_n+t_n.
\end{split}
\end{equation}
At the same time
\begin{equation}\notag
\begin{split}
\sigma_n(x(t))&=\min_{i_1,\dots,i_n} \dist(x(t),\spn\{g_{i_1},\dots,g_{i_n}\}) \\
&\geq\dist(x(t),\spn\{e_1,\dots,e_n\})-\varepsilon=d_n+t_n-\varepsilon,
\end{split}
\end{equation}
since the nearest point to $x(t)$ in each $\spn\{g_{i_1},\dots,g_{i_n}\}$ has the norm  at most $|x(t)|=d_0<M$.

According to  Lemma \ref{lemmaB} there is $x\in \R^{N+1}$ with the  $n$-term deviations with respect to $D$  equal to $\{d_n\}$.

\end{proof}

\section{Appendix: special bases}

In this section we construct the  bases with special geometrical properties we needed above to build examples of various dictionaries. The interested reader is invited
to come up with constructions of his own.

\begin{construction}\label{matrix}
For every $\eps>0$ and $N\in \NN$ there is a normalized basis $\{g_i\}_{i=1}^N$ of $\R^N$ so that for every $n\in\{1,\dots,N\}$
\begin{enumerate}
\item[{\rm(i)}]  $\spn\{g_1,\dots,g_n\}=\spn\{e_1,\dots,e_n\}$, which we denote by $\R^n$;
\item[{\rm(ii)}]    for every $n$-tuple of indices $i_1<\dots<i_n$ the subspace
$\spn\{g_{i_1},\dots,g_{i_n}\}$ is so close to $\R^n$ that the Hausdorff distance between their intersections with the  unit ball  is less than $\eps$.
\end{enumerate}
\end{construction}
\begin{proof}
Let $A$ be a lower triangular $(N\times N)$-matrix with the property that for every $n\in\{1,\dots,N\}$,  every $n\times n$ minor of $A$ formed by some $n$ of its rows and the first $n$ of its columns is not zero.
We claim that there are positive numbers $1=c_1\geq c_2\geq \dots \geq c_N>0$ so that if we for all $j$'s multiply the $j$-th column of $A$ by $c_j$ then the normalized rows of the so obtained  new matrix $\tilde A$ form the desired basis $\{g_i\}_{i=1}^N$.

First observe that matrices $A$ as above do exist. The family $\Aa$ of such matrices is in fact a full measure set in the space of all lower triangular $N\times N$ matrices:  every $n$-rows and  first-$n$-columns minor is a non-zero polynomial, hence it vanishes  on a set of measure zero.

Here are two  relevant properties of  every $A\in \Aa$.

If  we multiply each column of $A$ by a non-zero number then the new matrix is, clearly,  again in $\Aa$.
If we fix any $n$ rows of $A$ and restrict each of them to the first $n$ coordinates then the span of these new $n$ vectors is $\R^n$.

This property applied to the first $n$ rows together with $A$ being
lower triangular implies  the property (i) of the lemma.

To obtain the property (ii) we proceed by induction on $n$. First we multiply all columns except the first one by a non-zero positive number so small that (ii) is satisfied for $n=1$. Then we multiply all columns of the new matrix except for the first two columns by a non-zero positive number so small that (ii) is satisfied for $n=2$. And so on, till we after finitely many steps obtain the matrix $\tilde A$.
\end{proof}

The following elementary construction of bases almost contained in a hyperplane comes in handy  when building examples of dictionaries with $\rho=0$.
\begin{construction}\label{build}
Let $s$ be unit a vector in $\R^{m+1}$, and let $\eps>0$.
There exists a basis of $\R^{m+1}$ consisting   of $m+1$ unit vectors
$ g_0,\dots,g_m$
with the following property:  any  $m$-term linear combination $w$ of $g_0,\dots,g_m$ (i.e., $w=\sum_{j\not= i}\lambda_jg_j$ for some $i\in \{0,\dots,m\}$) satisfies the inequality
\begin{equation}\label{ws}
\left|\langle w, s\rangle\right|\leq \varepsilon|w|.
\end{equation}
More generally, let $\SSS=\{s_\nu\}_{\nu=1}^\infty$ be  a countable  set of  non-zero elements of $\ell_2$.  For each $\nu\in \NN$, let   $\eps(\nu)>0$, let
  $n(\nu)\in \NN$,   and
  let     $H_\nu$  be a $n(\nu)$-dimensional Euclidean space. Let
$H$ be the separable Hilbert space
$$
H=l_2\oplus_\perp H_1\oplus_\perp H_2\oplus_\perp\dots.
$$
In each subspace  $\spn\{s_\nu\}\oplus_\perp H_\nu$ of $H$ there is  a basis of $n(\nu)+1$ unit vectors $g_0^\nu,\dots,g_{n(\nu)}^\nu$
with the following property. Any  $m$-term linear combination $w$ of $g_j^\nu$'s (i.e., $w= \lambda_{j_1}g_{j_1}^{\nu_1}+\dots+\lambda_{j_m}g_{j_m}^{\nu_m}$)  with all the relevant
$n(\nu)\geq m$  and all the relevant $\eps(\nu)<\eps$ satisfies the inequality
$$
|Pw|\leq \varepsilon\sqrt m|(Id-P)w|;
$$
here $P:H\to \ell_2$ is the orthogonal projection.
\end{construction}
\begin{proof}
Assume $0<\eps<1/2$. We choose an orthonormal basis $g_1,\dots,g_m$  of $s^\perp$. The  vector $g_0$
will be only almost orthogonal to $s$, namely
$$
g_0=\left(\sqrt{1-\varepsilon^2}(g_1+\dots +g_m)+\varepsilon s\right)(m(1-\varepsilon^2)+\varepsilon^2)^{-1/2}.
$$
It is readily computed that
$\left|\langle w, s\rangle\right|\leq \varepsilon|w|$ and hence also
$\left|\langle w, s\rangle\right|\leq 2\varepsilon|h|$, where $h$ is the orthogonal projection of $w$ onto $s^\perp$.

For the second part of the construction, we can assume that $\SSS$ is contained in the unit sphere of $\ell_2$, and for each $\nu$ define the basis $g_0^\nu,\dots,g_{n(\nu)}^\nu$ of $\spn\{s_\nu\}\oplus_\perp H_\nu$ the same way as above, with $\varepsilon=\varepsilon(\nu)$ in (\ref{ws}).

Let $w= \lambda_{j_1}g_{j_1}^{\nu_1}+\dots+\lambda_{j_m}g_{j_m}^{\nu_m}$ be given. We group the terms with equal $\nu_j$, so that $w=\sum_{k=1}^{\leq m}  w_k$, where  each $w_k$ is an at most $m$-term linear combination of the just defined  basis of $\spn\{s_{\nu_k}\}\oplus_\perp H_{\nu_k}$; the    $\nu_k$'s are different,  all of them satisfy  $n(\nu_k)\geq m$, and  $\eps(\nu_k)<\eps$.
Then, since the vectors $(Id-P)w_k$ are pairwise orthogonal,
\begin{equation}\notag
\begin{split}
|Pw|&\leq \sum_{k=1}^{\leq m} |P w_k|=\sum_{k=1}^{\leq m}  |\langle w_k, s_{\nu_k}\rangle|\leq 2\eps \sum_{k=1}^{\leq m} |(Id-P)w_k| \\
&\leq 2\varepsilon\sqrt m|(Id-P)w|.
\end{split}
\end{equation}
\end{proof}

\subsection*{Acknowledgements}
We thank V.N. Temlyakov, S.V. Konyagin, A.P. Starovoitov and  A.A. Pekarskii for fruitful discussions.

\end{document}